\documentclass[11pt]{amsart}
\usepackage{amsfonts,eucal,amssymb}
\usepackage[all]{xy}
\begin{document}

\newcommand{\SmCorr}{{\sf SmCorr}}
\newcommand{\aA}{{\forall A}}
\newcommand{\Alg}{{\rm Alg}}
\newcommand{\cCorr}{{\sf Corr}}
\newcommand{\mE}{{\mathcal{E}}}
\newcommand{\mR}{{\mathcal{R}}}
\newcommand{\diff}{{\rm Diff}}
\newcommand{\bA}{{\bf A}}
\newcommand{\mtm}{{\sf{MTM}}}      
\newcommand{\odd}{{\rm odd}}
\newcommand{\Oh}{{\rm O}}
\newcommand{\Q}{{\mathbb Q}}
\newcommand{\R}{{\mathbb R}}
\newcommand{\Z}{{\mathbb Z}}
\newcommand{\alg}{{\rm alg}}
\newcommand{\C}{{\mathbb C}}
\newcommand{\Corr}{{\rm Corr}}
\newcommand{\Tor}{{\rm Tor}}
\newcommand{\Ext}{{\rm Ext}}
\newcommand{\Hom}{{\rm Hom}}
\newcommand{\Maps}{{\rm Maps}}
\newcommand{\Mod}{{\rm Mod}}
\newcommand{\MU}{{\rm MU}}
\newcommand{\Comod}{{\rm Comod}}
\newcommand{\reps}{{\rm reps}}
\newcommand{\ess}{{\mathbb S}}
\newcommand{\Spf}{{\rm Spec}}
\newcommand{\TC}{{\rm TC}}
\newcommand{\F}{{\mathbb F}}
\newcommand{\G}{{\mathbb G}}
\newcommand{\Wh}{{\rm Wh}}

\title {A theory of base motives}
\author{Jack Morava}
\address{Department of Mathematics, Johns Hopkins University, Baltimore,
Maryland 21218}
\email{jack@math.jhu.edu}
\thanks{This research was supported by DARPA and the NSF}
\date {20 August 2009}

\begin{abstract}
\noindent
When $A$ is a commutative local ring with residue field $k$, the derived tensor 
product
\[
- \otimes^L_A k : D(A) \to (k-\Mod)
\]
lifts to a functor taking values in a category of modules over the
`Tate cohomology' $R\Hom^*_A(k,k)$, which is the universal enveloping 
algebra of a certain Lie algebra. Under reasonable conditions this lift 
satisfies a spectral sequence of Adams (or Bockstein) type. \medskip

\noindent
In a suitable category of ring-spectra, replacing $A \to k$ by $A(*) \to \ess$ or 
$\TC(\ess) \to \ess$ yields interesting Hopf objects, with Lie algebras free after 
tensoring with $\Q$, analogous to those of motivic groups studied recently by Deligne,
Connes and Marcolli, and others. \end{abstract}

\maketitle

\noindent
[This is a sequel to and continuation of a talk at last summer's conference in Bonn
honoring Haynes Miller [23]. I owe many mathematicians thanks for helpful conversations and
encouragement, but want to single out John Rognes particularly, and thank him as well for
organizing this wonderful conference.] \bigskip

\begin{center}
{\bf \S 1 Prologue} 
\end{center} \bigskip

\noindent
Historically, the first part of the stable homotopy ring to be systematically understood
was the image of the $J$-homomorphism
\[
J : \pi_{k-1}\Oh = KO_k(*) \to \lim_{n\to \infty}\pi_{n+k-1}(S^n) = \pi^S_{k-1}(*) \;,
\]
defined on homotopy groups by the map
\[
\Oh_n \to \Oh \to \lim_{n \to \infty} \Omega^{n-1}S^{n-1} := Q(S^0) \;;
\] 
it factors through 
\[
KO_{4k} = \Z \to \zeta(1-2k) \cdot \Z/\Z \subset \pi^S_{4k-1}(*)  
\]
(at least, away from two). \bigskip

\noindent
In more geometric terms, a real vector bundle over $S^{4k}$ defines a stable cofiber sequence
\[
\xymatrix{
S^{4k-1} \ar[r]^\alpha & S^0 \ar[r] & {\rm cof} \; \alpha \ar[r] &  S^{4k} \cdots }
\]
and hence an extension
\[
[0 \to KO(S^{4k}) \to \cdots \to KO(S^0) \to 0] 
\]
in the group
\[
\Ext_{\rm Adams}(KO(S^0),KO(S^{4k}) \cong H^1_c(\hat{\Z}^\times,\hat{K}O(S^{4k})) \;,
\]
where the Adams operation $\psi^\alpha, \; \alpha \in \hat{\Z}^\times$ acts on
$\hat{K}O(S^{2k})$ by $\psi^\alpha(b^k) = \alpha^k b^k$. This (essentially Galois)
cohomology can be evaluated, via von Staudt's theorem, in terms of Bernoulli numbers. \bigskip
                                   
\noindent
In the arithmetic-geometric context, Deligne and Goncharov [[11]; cf [25] for a more
homotopy-theoretic account] have constructed an abelian tensor $\Q$-linear category $\mtm$ of
{\bf mixed Tate motives} over $\Z$, generated by objects $\Q(n)$ satisfying a (small, ie trivial
when $*>1$) Adams-style spectral sequence
\[
\Ext_\mtm^*(\Q(0),\Q(n)) \Rightarrow K_{2n-*}(\Z) \otimes \Q
\]
The groups on the right have rank one in degree $4k+1$, with generators corresponding (via
Borel regulators) to $\zeta(1+2k)$. \bigskip

\noindent
{\bf These same zeta-values} appear in differential topology [18] in the classification of smooth
(`Euclidean') cell bundles over the $4k+2$-sphere. There, both even and odd zeta-values can be
seen as having a common origin, summarized by a diagram (where, implicitly, $n \to \infty$)
\[
\xymatrix{
{} & X \ar@{.>}[dl] \ar[d] \ar@{-->}[dr] & {}\\
B\Oh_n \ar[d] \ar[r] & B\diff(E^n) \ar[r] & B\diff_c(\R^n) \ar[d] \\
BQ(S^0) & {} & \Omega \Wh(*) \;. }
\]

\noindent
The space $\Wh(*)$ on the bottom right is Waldhausen's smooth pseudoisotopy space, which
appears in 
\[
K(\ess) = A(*) = \ess \vee \Wh(*) \;.
\]
The shift by a double suspension in the cell versus vector bundle story is explained by the factor
$B$ on the lower left, and $\Omega$ on the lower right. The odd zeta-values appear in both
geometry and topology because the natural map 
\[
K(\Z) \to K(\ess)
\]
is a rational equivalence.

\newpage

\noindent
This suggests that some of the ideas of differential topology might be usefully reformulated 
in terms of a category of `motives over $\ess$' analogous to the arithmetic geometers' motives
over $\Z$, with the algebraic $K$-spectrum of the integers replaced by Waldhausen's
$A$-theory: these zeta-values might then provide a trail of breadcrumbs leading us to some
deeper insights. \bigskip

\noindent 
In the following section we recall some machinery from homological algebra, regarding 
\[
K(\ess) = A(*) = \ess \vee \Wh(*) \to \ess
\]
and 
\[
\TC(\ess) \sim \ess \vee \Sigma \C P^\infty_{-1} \to \ess
\]
(mod completions) as analogs of local rings over $\ess$, with the appropriate trace maps
interpreted as quotients by maximal ideals.  Note that the algebraic $K$-theory spectrum of $\Z$
lacks such an augmentation. \bigskip

\noindent
Tannakian formalism identifies the category $\mtm$ of mixed Tate motives as representations of
a certain pro-affine $\Q$-groupscheme with free graded Lie algebra, conjecturally related to other
areas of mathematics such as algebras of multiple zeta-values and renormalization theory [8, 10].
In the context proposed here, a similar group object 
\[
\Spf \; \ess \wedge_A \ess
\]
appears as derived automorphisms of $A$. $\S 3$ proposes to define a cycle map from
arithmetic motives to their $A$-theoretic analogs, conjecturally identifying these arithmetic and
geometric motivic groups. 
\bigskip \bigskip

\begin{center}
{\bf \S 2 Brave new local rings}
\end{center} \bigskip
                              
\noindent
{\bf 2.1} I'll start with work on commutative local rings, eg  $A \to k$ with maximal
ideal $I$, with roots in the very beginnings [27] of homological algebra. Eventually $A$ will be
graded, or a DGA. \bigskip

\noindent
The functors
\[
H_*(A,-) := \Tor_*^A(k,-)
\]
and 
\[
H^*(A,-) := \Ext^*_A(k,-) 
\]
appear in Cartan-Eilenberg; the first is covariant, and the second is contravariant, in $A$.
I'll be concerned mostly with 
\[
H_*(A,k) = \Tor^A_*(k,k) \; {\rm and}
\] 
\[
H^*(A,k) = \Ext^*_A(k,k)\;.
\]
Under reasonable finiteness conditions, these are dual $k$-vector spaces [the associativity sseqs
[7 XVI \S 4] degenerate]: in fact they are dual Hopf algebras, with $H^*(A,k)$ being the
universal enveloping algebra of a graded Lie algebra [2]. \bigskip

\noindent
{\bf 2.1.1 ex:} If $A = \Z_p \to \F_p$ is the residue map then 
\[
H_*(A,k) = E_*(Q_0)
\] 
is an exterior algebra on a Bockstein element of degree one. If $A = k[\epsilon]/(\epsilon^2)$
then $H_*(A,k) = k[x]$ ($|x|=2$) is the Hopf algebra of the additive group. These are the first
manifestations of Koszul duality.\bigskip

\noindent
{\bf 2.1.2 remarks:} For local rings this homology is closely related to Hochschild theory [7 X
\S 2], so it may also be related to recent work [3, 16] on Hopf algebra structures on THH.
\bigskip

\noindent
{\bf 2.1.3 Proposition:} The homological functor
\[
M \mapsto H_*(M \otimes^L_{\Z_p} \F_p) := \overline{M} : D(\Z_p-\Mod) \to \F_p-\Mod
\]
lifts to the category of $E(Q_0)$-comodules. There is a Bockstein spectral sequence
\[                  
\Ext^*_{E(Q_0)-\Comod}(\overline{M},\overline{N}) \Rightarrow \Hom_{D(\Z_p-
\Mod)}(M,N) \;.
\]
\bigskip

\noindent
{\bf 2.1.4 Definition:} $\G(A) := \Spf \; H_*(A,k)$ is an affine (super) $k$-groupscheme; its
grading is encoded by an action of the multiplicative group 
\[
\G_m = \Spf \; k[\beta^{\pm 1}] \;, 
\]
and $\tilde{\G}A) := \G(A) \rtimes \G_m$. \bigskip

\noindent
{\bf 2.2.1} The Bockstein spectral sequence generalizes: if $M \in D(A-\Mod)$, let
\[
\overline{M} = H_*(M \otimes^L_A k) = H_*(M \otimes_A \bA) \in (k-\Mod) \;,
\]
where $A \to \bA \to k$ is a factorization of the quotient map through a cofibration and a weak
equivalence (ie $\bA$ is a resolution of $k$, eg 
\[
\xymatrix{ 0 \ar[r] & \Z_p \ar[r]^p & \Z_p \ar[r] & 0 \; ) \;.}
\]  

\noindent
{\bf Proposition:} The functor $M \to \overline{M}$ lifts to a homological functor
\[
D(A-\Mod) \to (\tilde{G}(A)-\reps) \;,
\]                               
and there is an `ascent' sseq
\[
\Ext^*_{\tilde{G}(A)-\reps}(\overline{M},\overline{N}) \Rightarrow
\Hom_{D(A-\Mod)}(M,N)
\]
of Adams (alt: Bockstein) type \dots \bigskip

\noindent
The {\bf Proof} is as in Adams' Chicago notes [1], replacing the map $\ess \to \MU$ with $A \to
k$: thus $MU_*(X)$ becomes an $MU_*MU$-comodule by taking homotopy groups of the
composition
\[                                 
X \wedge \MU = X \wedge \ess \wedge \MU \to X \wedge \MU \wedge \MU 
\]
\[
= (X \wedge \MU) \wedge_\MU (\MU \wedge \MU) \;,
\]
yielding
\[
MU_*(X) \to MU_*(X) \otimes_{MU_*} MU_*MU \;.
\]
In the present context the comodule structure map comes from taking the homology of the
composition
\[
M \otimes_A \bA = M \otimes_A A \otimes_A \bA \to M \otimes_A \bA \otimes_A \bA 
\]
\[
= (M \otimes_A \bA) \otimes_\bA (\bA \otimes_A \bA) \;,
\]
resulting in $\overline{M} \to \overline{M} \otimes_kH_*(A,k) \; . \; \Box$ \bigskip

\noindent
{\bf 2.2.2 ex:} The bar construction provides a cofibrant replacement for $k$, with underlying
algebra
\[
\oplus_{n \geq 0} \otimes^n I[1]
\]
and a suitable differential. When $A = k \oplus I$ is a singular ($I^2 = 0$) extension, the
differential is trivial, and $\Ext^*_A(k,k)$ is the universal enveloping algebra of the free Lie
algebra on $I^*[1]$. \bigskip

\noindent
{\bf 2.2.3} {\bf Convergence} of such generalized Adams spectral sequences is a complicated
topic, related to extending the Tannakian formalism when there may be inequivalent fiber
functors [22]. The stable homotopy category is very unlike that of pure motives, which is
semisimple in interesting cases: instead, stable homotopy is more like the categories of $\F_p$-
representations of finite $p$-groups, whose structure is encoded entirely through iterated
extensions of trivial objects. Away from characteristic zero, it is often unrealistic to hope to
recover the full structure of an abelian (or triangulated) monoidal category in terms of the
automorphism group of a fiber functor; instead one usually gets at best a spectral sequence which
may allow the recovery of the graded object associated to a filtration of some localization of the
original category. \bigskip

\noindent
The generalized fiber functor defined by topological $K$-theory, for example, has ${\rm
Gal}(\Q_{\rm ab}/\Q) \cong \hat{\Z}^\times$ as (more or less) its motivic group, and the
associated spectral sequence `sees' only the image of the $J$-homomorphism; other fiber
functors see different parts of (some generalization [17] of) the prime ideal spectrum of the stable
homotopy category. One of the more interesting issues emerging from this picture is the relation
of deformations of fiber functors (eg, taking values in categories of modules over a local ring)
and their motivic groups. \bigskip

\noindent
{\bf 2.3.1} Examples closer to homotopy theory appear in recent work of Dwyer, Greenlees, and
Iyengar. Suppose for example that $X$ is a connected pointed space (eg with finitely many cells
in each dimension), and let $X_+ = X \vee S^0$ be $X$ with a disjoint basepoint appended. Its
Spanier-Whitehead dual 
\[
X^D := \Maps_\ess(\Sigma^\infty X_+,\ess) 
\]
is an $E_\infty$ ring-spectrum, with augmentation $X^D \to \ess$ given by the basepoint.
\bigskip

\noindent
The Rothenberg-Steenrod construction [14 \S 4.22] then yields an equivalence
\[
\Hom_{X^D-\Mod}(\ess,\ess) \sim \ess[\Omega X]
\]
of ($A_\infty, \; {\rm co}-E_\infty)$ Hopf algebra objects in the category of spectra. \bigskip

\noindent
If $X$ is simply connected, there is a dual result with coefficients in the Eilenberg-MacLane
spectrum $k = Hk$ of a field: then the `double commutator'
\[
\Hom_{k[\Omega X]-\Mod}(k,k) \sim C^*(X,k)
\]
is homotopy equivalent to the (commutative) cochain algebra of $X$. [This puts the homotopy
groups
\[
\pi_*C^*(X,k) \cong H^{-*}(X,k)
\]
in negative dimension.] This sharpens a classical [4] analogy between the homology of
loopspaces and local rings. \bigskip                        

\noindent
The (generalized Koszul duality?) functor 
\[ 
M \mapsto \Hom_{X^D}(M,\ess) : (X^D-\Mod) \to (\ess[\Omega X]-\Mod) \;.
\]
seems worth further investigation \dots \bigskip

\noindent
{\bf 2.3.2 ex:} Suspensions are formal, so if $k = \Q$ and $X = \Sigma Y$ then 
\[
X^D \otimes \Q \sim H^{-*}(\Sigma Y,\Q)
\]
is a singular extension of $\Q$, so $\Q[\Omega \Sigma Y]$ is the universal enveloping algebra of
the free Lie algebra on the graded dual of $\tilde{H}^{-*-1}(Y,\Q)[1]$. \bigskip

\noindent
Recent work of Baker and Richter [5] identifies the Hopf algebra of noncommutative symmetric
functions with the integral homology $H_*(\Omega \Sigma \C P^\infty)$ as the universal
enveloping algebra of a free graded Lie algebra. The dual Hopf algebra $H^*(\Omega \Sigma \C
P^\infty)$ is the (commutative) algebra of quasi-symmetric functions.  
\bigskip

\noindent
{\bf 2.3.3} The topological cyclic homology $\TC(\ess;p)$ of the sphere spectrum (at $p$) is an
$E_\infty$ ringspectrum, equivalent to the $p$-completion of $\ess \vee \Sigma \C
P^\infty_{-1}$ [20]; the subscript signifies a twisted desuspension of projective space by
the Hopf line bundle. \bigskip

\noindent
From now on I'll be working over the rationals, eg with the graded algebra
\[
\TC_{2n-1}(\ess;\Q_p) \cong \Q_p \oplus \Q_p \langle e_{2n-1}\rangle \;,
\]
$n \geq 0$ (with trivial multiplication). \bigskip

\noindent
{\bf 2.4.1} The multiplication on a ring-spectrum $A$ defines a composition
\[
[X,A \wedge Y] \wedge [Y, A \wedge Z] \to [X, A \wedge Z]
\]
(on morphism objects in spectra) by
\[ 
X \to A \wedge Y \to A \wedge A \wedge Z \to A \wedge Z \;.
\]
The map $X \to A \wedge X$ defines a functor from the category with $\ess$-modules (eg
$X,Y$) as objects, and 
\[
\Corr_A(X,Y) := [X,A \wedge Y]                                   
\]
as morphisms, to the category of $A$-modules, because
\[
\Corr_A(X,Y) = [X,A \wedge Y] \to [X,[A,A \wedge Y]_A]
\]
\[
\cong [A \wedge X, A \wedge Y]_A \;.
\]

\noindent
Let $(A-\Corr)$ be the triangulated subcategory of $(A-\Mod)$ generated by the image of this
construction. The augmentation of $A$ defines a functor from $(A-\Corr)$ to $\ess$-modules
which is the identity on objects, and is given on morphisms by
\[
[X,A \wedge Y] \to [X,\ess \wedge Y] = [X,Y] \;.
\]
I propose to inherit the composition of this functor with rationalization as an analog of the `fiber
functor' in \S 2.1: \bigskip

\noindent
{\bf Corollary:} This homological functor $(A-\Mod) \to (\Q-\Mod)$ lifts to the category of 
$\tilde{\G}(A \otimes \Q)$-representations, yielding a spectral sequence
\[
\Ext^{*,*}_{\tilde{\G}(A \otimes \Q)-\reps}(X,Y) \Rightarrow \Corr^*_A(X,Y) \;.
\] 

\noindent
{\bf Proof:} $(X \wedge A) \wedge_A \ess = X \dots \; \Box$ \bigskip

\noindent
{\bf 2.4.2} A free Lie algebra has cohomological dimension one, so when $A$ is $\TC(\ess;p)$
and $X$ and $Y$ are spheres, this spectral sequence degenerates to 
\[
\Ext^1_ {\tilde{\G}(\TC \otimes \Q_p)}(S^{2n}_{\Q_p},S^0_{\Q_p}) \cong
\TC_{2n-1}(\ess,\Q_p) \;,
\]
with left-hand side isomorphic to 
\[
H^{1,0}_{\rm Lie}(\mathfrak{F}(\widetilde{\TC}^*[1]),S^{-2n}_{\Q_p}) \cong 
\Hom^0(\widetilde{\TC}^*[1],S^{-2n}_{\Q_p}) \;, 
\]
which is just the one-dimensional vector space
\[
\Q_p \langle e_{2n-1}[1]b^{-n} \rangle \;.
\]

\noindent
At a regular odd prime $p$ (cf. [13, 24]), 
\[
\Wh(*)/\Sigma {\rm coker} \; J \; \sim \; \Sigma {\mathbb H} P^\infty \;;
\]
the cokernel of the $J$-homomorphism is a torsion space, so this yields a spectral sequence 
\[
\Ext^*_{\tilde{\G}(A \otimes \Q)}(S^{2n}_\Q,S^0_\Q) \Rightarrow A_{2n-*}(*) \otimes \Q 
\]
with $A_{4k+1}(*) \otimes \Q \cong \Q \langle e_{4k+1}[1]b^{-2k-1} \rangle$ .

\newpage

\begin{center}          
{\bf \S 3 $A$-theoretic motives}
\end{center} \bigskip

\noindent
{\bf 3.1} A retractive space $Z$ over $X$ is a diagram 
$\xymatrix{X & Z \ar[l]_r & X \ar[l]_s }$
which composes to the identity $1_X$: it's a space over $X$ with a cofibration section. $Z$ is
said to be finitely dominated if some finite complex is retractive over it [15]. \bigskip 

\noindent
Waldhausen showed that finitely dominated retractive spaces over $X$ form a category with
weak equivalences and cofibrations, and that the $K$-theory spectrum $A(X)$ of this category
can be identified with $K(\ess[\Omega X])$. \bigskip

\noindent
More generally, $Z$ is relatively retractive over $X$, with respect to a map $p:X \to Y$, if the
homotopy fiber of $p \circ r$ over any $y \in Y$ is finitely dominated as a retractive space over
the homotopy fiber of $p$ above $y$. The category $\mR(p)$ of such spaces is again closed
under cofibrations and weak equivalences, with an associated $K$-theory spectrum $A(X \to
Y)$. \bigskip

\noindent
Bruce Williams [28 \S 4] (using a formalism developed in algebraic geometry by Fulton and
MacPherson) shows that this functor has a rich bivariant structure: compositions 
\[
A(X \to Y) \wedge A(Y \to Z) \to A(X \to Z) \;,
\]
good behavior under products, \&c. It behaves especially well on fibrations; in particular, the
spectra 
\[
\aA(X,Y) := A(X \times Y \to X) 
\]
(defined by relatively retractive spaces $Z$ over $X\times Y \to X$) admit good products 
\[
\aA(X,Y) \wedge \aA(Y,Z) \to \aA(X,Z) \;.
\]

\noindent
Let $A-\cCorr$ be the triangulated envelope [6] of the symmetric monoidal additive
category with finite CW complexes $X,Y$ as objects, and $\aA_0(X,Y) = \pi_0\aA(X,Y)$ as
morphisms. Composition
\[
A(X \times Y \to Y) \to [X,A(Y)] \to [X,A \wedge Y]
\]
of the standard assembly map with a slightly less familiar relative co-assembly map [12 \S 5]
defines a monoidal stabilization functor 
\[
(A-\cCorr) \to (A-\Corr)
\]
analogous to inverting the Tate motive, or to the introduction of desuspension in classical
homotopy theory. However, $A$-theory of spaces is a highly nonlinear functor, and might
possess other interesting stabilizations. \bigskip

\noindent
{\bf 3.2} The motivic constructions of Suslin and Voevodsky [27] begin with a category whose
objects are schemes of finite type over some nice base, and whose morphism groups
$\SmCorr(V,W)$ of (roughly) sums of irreducible subvarieties $Z$ of $V \times W$ which are
finite with respect to the projection $V \times W \to V$, and surjective on components of $V$.
\bigskip

\noindent
When $V$ and $W$ are defined over a number field (eg $\Q$), classical arguments [cf. eg [21]]
show that 
\[
Z(\C) \cup V(\C) \times W(\C) \to V(\C) \times W(\C)
\]
is finitely dominated relatively retractive with respect to $V(\C) \times W(\C) \to V(\C)$, 
defining a cycle class homomorphism
\[
\SmCorr(V,W) \to \aA_0(V(\C),W(\C)) \;,
\]
and hence a functor 
\[
V \mapsto V(\C) : (\SmCorr) \to (A-\cCorr) \;.
\] 

\noindent
My {\bf hope} is that this will lead to an identification of the motivic group for the category of
mixed Tate motives with $\tilde{\G}(A \otimes \Q)$. It seems at least possible that 
$\tilde{\G}(\TC \otimes \Q)$ is the larger motivic group seen in physics 
[8, 9 \S 3.1] by Connes and Marcolli. \bigskip

\noindent
{\bf 3.3} I don't want to end this sketch without mentioning one last possibility. Dundas and
$\O$stv$\ae$r have proposed a bivariant K-theory based on categories $\mE(E,F)$ of suitably
exact functors between the categories of (cell) modules over (associative) ring-spectra $E$ and
$F$. \bigskip

\noindent
These module categories are to be understood as categories with weak equivalences and
cofibrations; the exact functors are to preserve these structures, and be additive in a certain sense.
$\mE (E,F)$ is again a Waldhausen category, which suggests that the category $(\Alg_A)$
with associative ring-spectra $E,F$ as its objects, and 
\[
\Alg_A(E,F) := K(\mE (E,F))
\]
as morphisms, is an interesting analog of categories of noncommutative correspondences
proposed by various research groups [9 \S 6, 19 \S 4]. It seems reasonable to expect
that this category will naturally be be enriched over $A$. \bigskip

\noindent
A space $W$ over $X \times Y$ defines an $X^D$-$Y^D$ bimodule $W^D$, and 
\[
W \mapsto \Hom_{X^D-\Mod}(W^D,-)
\]
is a natural candidate for an exact functor, and hence a map
\[
\mR(X \times Y \to X)  \rightarrow \mE(X^D,Y^D) \;.
\]
If so, this might define another interesting stabilization of $A-\cCorr$, related more closely to the
Waldhausen $K$-theory of Spanier-Whitehead duals than to spherical group rings. \bigskip

\bibliographystyle{amsplain}

\end{document}